\documentclass[12pt,english]{amsart}
\usepackage[T1]{fontenc}
\usepackage[latin9]{inputenc}
\usepackage[letterpaper]{geometry}
\usepackage{amsthm}
\usepackage{amsmath}
\usepackage{amssymb}
\usepackage[unicode=true]
 {hyperref}
\usepackage{breakurl}

\makeatletter
\numberwithin{equation}{section}
\numberwithin{figure}{section}
\theoremstyle{plain}
\newtheorem{thm}{Theorem}[section]
  \theoremstyle{definition}
  \newtheorem{defn}[thm]{Definition}
  \theoremstyle{remark}
  \newtheorem{rem}[thm]{Remark}
  \theoremstyle{plain}
  
  \theoremstyle{plain}
  
  \theoremstyle{plain}
  
 \theoremstyle{definition}
  \newtheorem{example}[thm]{Example}

\usepackage{color}

\usepackage{babel}

\makeatother

\usepackage{babel}
\begin{document}

\title{Nigel Kalton and the interpolation theory of commutators}

\author{{Michael Cwikel}}

\address{Cwikel: Department of Mathematics, Technion - Israel Institute of
Technology, Haifa 32000, Israel }

\email{mcwikel@math.technion.ac.il }

\author{{Mario Milman}}

\address{Milman: Department of Mathematics, Florida Atlantic University, Boca
Raton, FL 33431, USA }

\email{mario.milman@gmail.com}

\author{{Richard Rochberg}}

\address{Rochberg: Department of Mathematics, Washington University, St. Louis,
MO, 63130, USA }

\email{rr@math.wustl.edu}

\thanks{The first author's work was supported by the Technion V.P.R.\ Fund
and by the Fund for Promotion of Research at the Technion. The second
author's work was partially supported by a grant from the Simons Foundation
(\#207929 to Mario Milman). The third author's work was supported
by the National Science Foundation under Grant No. 1001488.}

\begin{abstract}
{This is the second of a series of papers surveying
some small part of the remarkable work of our friend and colleague
Nigel Kalton. We have written it as part of a tribute to his memory.
It does not contain new results. One of the many topics in which Nigel
made very significant and profound contributions deals with commutators
in interpolation theory. It was our great privilege to work with him
on one of his many papers about this topic. Our main purpose here
is to offer} an introduction to that paper: A unified theory of commutator
estimates for a class of interpolation methods. Adv. Math. 169 (2002),
no. 2, 241--312. We sketch the theory of interpolation spaces constructed
using pseudolattices which was developed in that paper and which enables
quite general formulation of commutator theorems. We seek to place
the results of that paper in the general context of preceding and
subsequent
research on this topic, also indicating some applications to other
fields of analysis and possible directions for future research.
\end{abstract}
\maketitle

\section{Introduction}

Given an initial pair of spaces, the classical interpolation theories
provide us with general methods to construct new spaces that have
the \textit{interpolation property}, meaning that operators which
are bounded on both the initial spaces are also bounded on the newly
constructed spaces. A salient feature of these constructions is their
generality: They apply equally well to all the operators that satisfy
the prescribed initial assumptions. However, certain important operators
in Analysis \textbf{(}some of which are nonlinear\textbf{)} can be
shown to be bounded only because of subtle cancellations. For many
of these operators the direct application of standard interpolation
methods turned out to be inadequate, and the challenge was how to
extend the classical theory to cope with this situation.

Although much more still needs to be done, in the last 30 years, through
the work of a number of mathematicians, among whom Nigel Kalton was
one of the leaders, we now have an extensive body of knowledge on
how to handle commutators, including certain nonlinear ones. Moreover,
the theory has naturally expanded to provide a general framework to
study certain types of cancellations. Many new applications in disparate
fields have been developed using these ideas.

Informally, the operators that originated the interpolation theory
of commutators are of the form: 
\[
\left[T,\Omega\right]=T\Omega-\Omega T,
\]
 where $T$ is a bounded linear operator (a ``good\textquotedblright{}\ operator)
in some interpolation scale and $\Omega$ is a possibly nonlinear
operator which is unbounded on relevant space(s) of the interpolation
scale (a ``bad\textquotedblright{}\ operator). At first sight it
would seem that forming such a commutator is going to worsen our predicament,
since it will require us to cope with two unbounded operators instead
of one...but, somehow, the nasty behaviors of these two bad terms
cancel each other out. This saves the day, and provides a bounded
operator. Informally, we have the following scheme. 
\[
Bad-Bad=Good!
\]

As already hinted above, phenomena of cancellation can be quite subtle
and difficult to analyze. But it turns out to be possible, in most
of the instances that we will describe here, to handle them via appropriately
interpreted versions of the following simple fact: If two holomorphic
functions have a simple pole at the same point and their residues
at that point are equal, then their difference is holomorphic at that
point.

In the next section of this document we will deal mostly with the
background and some of the motivations for the theory presented in
\cite{ckmr}. There is considerable overlap with the introductory
section of \cite{ckmr} itself. Here we have allowed ourselves the
luxury of presenting some of the examples in more detail, and adding
a few more comments, some of them relating to some more recent developments.

Subsequently, in Section \ref{sec:UnifiedMethod}, we shall attempt
to give the flavor of the results of \cite{ckmr}. Together with Nigel,
we wrote \cite{ckmr} with the intention of making it applicable to
as wide a range of settings as possible. This means that the reader
must wade through several levels of definitions and conditions and
extra remarks about better constants and weaker hypotheses. We hope
that the simplified overview given in Section \ref{sec:UnifiedMethod}
will provide perspectives to make such a task somewhat easier.

Finally, in Section \ref{sec:FurtherDirections}, we will suggest
some possible directions for future research.

\section{Interpolation theory of commutators}

While the main focus of interpolation theory is the systematic construction
of spaces ``with the interpolation property\textquotedbl{} from given
ones (``interpolation methods\textquotedbl{}), it was observed in
\cite{rw} that associated with these methods there are certain (unbounded
and possibly nonlinear) operators $\Omega$ whose commutators with
bounded operators can be controlled inside the corresponding interpolation
scales. In due course it was realized that the resulting theory could
be expanded to provide a general framework to analyze and control
the underlying cancellations in different contexts.

\subsection{Commutators in the setting of the complex method of interpolation}
\begin{example}
\label{ex:MbBMO} (The Coifman-Rochberg-Weiss Commutator Theorem \cite{crw})
It is well known that the space $BMO=BMO(\mathbb{R}^{n})$ contains
unbounded functions, in particular, if $b$ is an element of $BMO\backslash L^{\infty},$
the multiplier operator $M_{b}(f)=bf$ is not bounded on, say, $L^{p}=L^{p}(\mathbb{R}^{n}).$
On the other hand, Calderón-Zygmund singular integral operators are
known to be bounded on $L^{p},$ $1<p<\infty.$ The commutator theorem
of Coifman-Rochberg-Weiss \cite{crw} asserts that the commutator
between a Calderón-Zygmund operator $T,$ and multiplication by a
$BMO$ function $b,$ is a bounded operator on $L^{p}$: 
\begin{equation}
\left[T,M_{b}\right]:L^{p}\rightarrow L^{p}.\label{commthm}
\end{equation}

\end{example}
In \cite{crw} one of the ways in which the authors proved the boundedness
of the commutator $[T,M_{b}]$ on $L^{p}$ was by identifying it (cf.~\cite[pp.~620--621]{crw})
as the derivative of a (Banach space valued) analytic function. They
also used some consequences of a connection (formulated below in Example
\ref{ex:bmo}) of $BMO$ with the $A_{p}$ classes of weights. 
\begin{rem}
\label{rem:ApplOne}We refer to \cite[pp.~257--258]{CLMS} for the
first of several applications in that paper of (\ref{commthm}). It
yields results about the properties of Jacobians of functions with
prescribed smoothness.\textbf{ }Related material can also be found
in Chapter 13 of \cite{IwaniecTMartinG2001}.
\end{rem}

In \cite{rw}, Rochberg and Weiss took the ideas of \cite{crw} one
step further and, in the process, initiated the general interpolation
theory of commutators. In particular, they introduced a class of operators
which act on (and are associated with the construction of) individual
complex interpolation spaces (cf.~\cite{ca}). Although operators
$\Omega$ in this class can be nonlinear and/or unbounded on their
designated complex interpolation space, their commutators $\left[T,\Omega\right]$
are always bounded operators on that same interpolation space, for
every choice of linear operator $T$ which acts boundedly on both
of the ``endpoint'' spaces. 

{The reader who is not yet familiar with Alberto
Calder\'on's complex method of interpolation can of course refer
to \cite{ca} or, for example, to Chapter 4 of \cite{bl}, or can
perhaps make do with a brief overview, such as the one given in our
earlier paper \cite{CwikelMMilmanMRochbergR2014} in this series.
}
Given a Banach pair $\vec{A}=(A_{0},A_{1})$, we consider the scale of complex
interpolation spaces $A_{t}=[A_{0},A_{1}]_{t}$, $t\in(0,1).$ We
fix a value of $t$ and, for each $x\in A_{t}$, let $f_{x,t}(z)$
be its representing function, i.e.\ $f_{x,t}\in\mathcal{F}(\vec{A})$
is such that%
\footnote{As in {\cite{CwikelMMilmanMRochbergR2014},} 
in this discussion we have chosen to ignore issues concerning
the existence and uniqueness of \textbf{$f_{x,t}$}. 
In the literature this difficulty is sometimes dealt with by working modulo bounded
operators.%
} 
\[
f_{x,t}(t)=x,\text{ and }\left\Vert f_{x,t}\right\Vert _{\mathcal{F}}=\left\Vert x\right\Vert _{A_{t}}.
\]
 Then we let $\Omega_{\vec{A}}:A_{t}\rightarrow A_{0}+A_{1},$ be
defined by 
\begin{equation}
\Omega_{\vec{A}}x=f_{x,t}^{\prime}(t).\label{eq:DefOmega}
\end{equation}
 The abstract interpolation theorem of Rochberg-Weiss now reads: if
$\vec{A}$ and $\vec{B}$ are Banach pairs and $T:\vec{A}\rightarrow\vec{B}$
is any bounded operator acting between them$,$ then for each $t\in(0,1)$
there exists a constant $c(t)>0$ such that 
\begin{equation}
\left\Vert \left[T,\Omega\right]x\right\Vert _{B_{t}}\leq c(t)\left\Vert T\right\Vert _{\vec{A}\to\vec{B}}\left\Vert x\right\Vert _{A_{t}},\label{cuatro}
\end{equation}
 where%
\footnote{During a first reading, you might prefer to just consider the special
case where $\vec{A}=\vec{B}$ to avoid the slight inconvenience of
having to deal with two different versions of $\Omega$. In fact,
via a straightforward use of direct sums, the general case can be
deduced from this special case.%
} $\left[T,\Omega\right]=T\Omega_{\vec{A}}-\Omega_{\vec{B}}T$. (One
should bear in mind that $\Omega_{\vec{A}}$ and $\Omega_{\vec{B}}$
depend on our choice of $t$.)

The idea of the proof is to exploit the cancellation that is exhibited
by the commutator $\left[T,\Omega\right]$. Given $x\in A_{t},$ let
$G(z)=(T(f_{x,t}(z))-f_{Tx,t}(z))/(z-t)$, then the cancellation $Tf_{x,t}(t)-f_{Tx,t}(t)=Tx-Tx=0,$
makes this function analytic also at $z=t$ with 
\begin{equation}
G(t)=\lim_{z\rightarrow t}G(z)=\lim_{z\rightarrow t}\frac{T(f_{x,t}(z))-f_{Tx,t}(z)}{z-t}=T\Omega_{\vec{A}}x-\Omega_{\vec{B}}Tx.\label{eq:vjta}
\end{equation}

The desired result (\ref{cuatro}) will now follow immediately, once
we show that $G\in\mathcal{F}(\vec{B})$. This is rather obvious,
but for later purposes we will be more explicit: We have obtained
$G$ from a function in $\mathcal{F}(\vec{B})$ by multiplying that
function by the scalar function $1/(z-t)$. We have seen that the
singularity of that scalar function at $z=t$ does not create any
problems. So the result follows trivially from the facts that $1/\left(z-t\right)$
is analytic everywhere else in the strip $\left\{ z\in\mathbb{C}:0<\mathrm{Re}\, z<1\right\} $
and is continuous and bounded on its closure.

\begin{example} \textbf{\label{ex:Weights}} Let $p\in[1,\infty)$
and let $w_{0}$ and $w_{1}$ be positive functions on $\mathbb{R}^{n}$.
Consider the pair of weighted $L^{p}$ spaces%
\footnote{In our discussion here, $L^{p}(w)$ on some measure space $\left(X,\Sigma,\mu\right)$
(usually $\mathbb{R}^{n}$ with Lebesgue measure) is defined by requiring
the norm $\left\Vert f\right\Vert =\left(\int_{X}\left|f\right|^{p}wd\mu\right)^{1/p}$
to be finite. But it should be noted that there are many papers in
which it is found to be more convenient to replace $w$ by $w^{p}$
in this definition. %
}$\vec{L^{p}}=(L^{p}(w_{0}),L^{p}(w_{1}))$. It is well known and easy
to show that the complex method of interpolation gives 
\begin{equation}
\left[L^{p}(w_{0}),L^{p}(w_{1})\right]_{t}=L^{p}(w_{0}^{1-t},w_{1}^{t})\mbox{ for each }t\in(0,1).\label{eq:com}
\end{equation}
 A simple calculation, which is in fact also a main step of the proof
of (\ref{eq:com}), shows that, for each $t\in(0,1)$ and for all
elements $x$ in a dense subspace of $\left[L^{p}(w_{0}),L^{p}(w_{1})\right]_{t}$,
we have $f_{x,t}(z)=(w_{1}/w_{0})^{(z-t)/p}x$. So, in this case,
for these elements%
\footnote{There are also other examples where it can happen that some general
explicit formula obtained for $\Omega x$ can be seen to coincide
with $f_{z,t}^{\prime}$ only for $x$ in some dense subspace of the
relevant space $A_{t}$. However this density may enable (\ref{cuatro}),
with $\Omega$ given by that explicit formula, to be deduced for all
$x\in A_{t}$. We are glossing over this issue here. %
}, the operator $\Omega_{\vec{L^{p}}}$ (cf.~(\ref{eq:DefOmega}))
is given by 
\begin{equation}
\Omega_{\vec{L^{p}}}=\frac{1}{p}x\log\left(w_{1}/w_{0}\right).\label{cinco}
\end{equation}
 (Various versions of this calculation, in this and other more elaborate
similar settings can be found in \cite{rw} pp.\ 335--340.)

\label{ex:bmo} This is an important special case of Example \ref{ex:Weights}.
If we now require $p>1$ and if $b$ is a $BMO$ function, then, as
more or less informally observed and/or used in a number of papers,
including \cite{CoiFef}, \cite{crw}, \cite{HuntMuckenhouptWheeden}
and \cite{rw}, etc.%
\footnote{In the case of functions of one variable this fact is linked to a
much earlier theorem of Helson-Szeg\H{o}. The authors of \cite{CoiFef}
also acknowledge contributions of R.~Gundy and of the authors of
\cite{HuntMuckenhouptWheeden}.%
}, and explained in detail, e.g., in \cite{GarsiaRubio} p.~409, there
exists $0<\alpha=\mathrm{O}(1/\left\Vert b\right\Vert _{BMO})$ (depending
on $p$) such that $e^{-\alpha b}$ and $e^{\alpha b}$ both belong
to the $A_{p}$ class of Muckenhoupt weights. Therefore, if $T$ is
a Calderón-Zygmund singular integral operator, then $T:L^{p}(e^{-\alpha b})\rightarrow L^{p}(e^{-\alpha b})$
and $T:L^{p}(e^{\alpha b})\rightarrow L^{p}(e^{\alpha b})$ continuously.
We have (cf.~(\ref{eq:com})) that $[L(e^{-\alpha b}),L^{p}(e^{\alpha b})]_{1/2}=L^{p}$
and so, by (\ref{cinco}), 
\[
\Omega f=\frac{1}{p}f\log\left(e^{\alpha b}/e^{-\alpha b}\right)=\frac{2}{p}f\alpha b=\frac{2}{p}\alpha M_{b}(f).
\]
 Consequently, (\ref{cuatro}) with $\vec{A}=\vec{B}=\vec{L^{p}}$
gives

\textbf{
\[
\left\Vert \left[T,M_{p}\right]f\right\Vert _{L^{p}}=\frac{p}{2\alpha}\left\Vert \left[T,\Omega\right]f\right\Vert _{L^{p}}\le c\left\Vert b\right\Vert _{BMO}\left\Vert f\right\Vert _{L^{p}},
\]
 }where the constant $c$ depends on $p$ and also on the norms of
$T$ as an operator on $L^{p}(e^{-\alpha b})$ and $L^{p}(e^{\alpha b})$.
Thus we recover the commutator theorem of \cite{crw}. \end{example}

\begin{example} \label{ex:L1Linfty}For the pair $\left(L^{1},L^{\infty}\right)$
we have $\left[L^{1},L^{\infty}\right]_{t}=L^{p}$ where $p=1/(1-t)$
for each $t\in(0,1)$. Then, for each simple function $x\in\left[L^{1},L^{\infty}\right]_{t}$
we have 
\[
f_{x,t}(z)=\left\Vert x\right\Vert _{p}{\displaystyle \frac{x}{\left|x\right|}\left(\frac{\left|x\right|}{\left\Vert x\right\Vert _{p}}\right)^{p(1-z)}}
\]
 and, consequently, 
\[
\Omega x=-px\log\left(\frac{\left|x\right|}{\left\Vert x\right\Vert _{p}}\right).
\]
 See Proposition 1.5 on p.~166 of \cite{MingFAN} for a detailed
discussion of a non-commutative version of this example. Cf.~also
a related calculation, with some applications, on pp.~315--318 of
\cite{rw}. \end{example}

\subsection{\label{sub:transops}Translation operators and commutators, and their
application to partial differential equations}

When studying operators on a scale of spaces it is also of interest
to consider ``translation'' operators and their corresponding commutators.
It should be made clear at the outset that here we are \textit{not}
talking about translation operators on function spaces, i.e., not
about maps of the kind $f\left(\cdot\right)\mapsto f(\cdot-x_{0})$
for some fixed $x_{0}$ in the underlying set. Instead it is the parameter
of interpolation $t$ which is being translated.

For the particular case of $L^{p}$ spaces this idea was developed
by Iwaniec \cite{i0} and Iwaniec and Sbordone \cite{is} (cf.~Example
\ref{iwansb} below) who found remarkable applications to nonlinear
PDEs (cf.~\cite{i0}, \cite{i}, and the references therein). More
generally, for the complex method of interpolation, these translation
operators can be defined as follows. Suppose that $x\in A_{t}$ is
represented by 
\[
x=f_{x,t}(t),\text{ with }\left\Vert x\right\Vert _{A_{t}}=\left\Vert f_{x,t}\right\Vert _{\mathcal{F}}.
\]
 Then, for sufficiently small $\varepsilon\in\mathbb{R}$, we define
\[
\mathcal{R}_{\varepsilon}x=f_{x,t}(t+\varepsilon).
\]
 It follows that 
\[
\mathcal{R}_{\varepsilon}:A_{t}\rightarrow A_{t+\varepsilon}.
\]
 The general version of the commutator theorem for translations (cf.~\cite{ckmr})
asserts that if $T:\vec{A}\rightarrow\vec{B}$ is a bounded linear
operator, then for all $t\in(0,1),$ there exists a constant $c=c(t)>0$
such that%
\footnote{Here, analogously to what happens in (\ref{cuatro}), $\mathcal{R}_{\varepsilon}$
has two different meanings in (\ref{iwan}), depending on whether
it acts on $A_{t}$ or on $B_{t}$. Here again one can readily deduce
the general case of this result from the special and more easily formulated
case where $\vec{A}=\vec{B}$. %
} 
\begin{equation}
\left\Vert \mathcal{R}_{\varepsilon}Tx-T\mathcal{R}_{\varepsilon}x\right\Vert _{B_{t+\varepsilon}}\leq c\left\vert \varepsilon\right\vert \left\Vert T\right\Vert _{\vec{A}\to\vec{B}}\left\Vert x\right\Vert _{A_{t}},\label{iwan}
\end{equation}
 where the constant $c$ depends only on $t$.

Of course, for each $x\in A_{t}$, we have $\lim_{\varepsilon\to0}(\mathcal{R}_{\varepsilon}x-x)/\varepsilon=\Omega x$
in suitable topologies (such as the norm topology in $A_{0}+A_{1}$)
which can lead us to think of (\ref{cuatro}) as a sort of limiting
case%
\footnote{This idea is brought into play in a rather more general setting in
Section 4 of \cite{krmi}.%
} of (\ref{iwan}). In fact quite generalized versions of (\ref{cuatro})
and (\ref{iwan}) (both formulated in Theorem 3.8 on page 256 of \cite{ckmr})
emerge easily and simultaneously from the same simple argument appearing
in the lower part of page 259 of \cite{ckmr}, once the required machinery
has been set up. \begin{example} \label{iwansb}(A commutator theorem
of Iwaniec-Sbordone \cite{is}) We consider ``translation operators''
in the setting of $L^{p}$ spaces. Let $1\leq p_{0}<p_{1}\leq\infty,$
be fixed, and let $\theta\in(0,1).$ Let $p\in(p_{0},p_{1})$ be given
by $1/p=\left(1-\theta\right)/p_{0}+\theta/p_{1}$ so that $L^{p}=[L^{p_{0}},L^{p_{1}}]_{\theta}.$
Moreover, for $z$ on the strip $0\le\mathrm{Re}\, z\le1$, we let
$h(z)=(1-z)/p_{0}+z/p_{1}$. Then, if $f$ is a non-zero element of
$L^{p},$ it can be represented in an optimal way, for the computation
of its norm in the interpolation space $L^{p},$ by the analytic function
\[
F_{f,\theta}(z)=\left\Vert f\right\Vert _{p}\left(\frac{\left\vert f\right\vert }{\left\Vert f\right\Vert _{p}}\right)^{h(z)p}\frac{f}{\left\vert f\right\vert },\text{ with }0\leq\mathrm{Re}\, z\leq1.
\]
 (The formula given above in Example \ref{ex:L1Linfty} is a limiting
case of this formula.) It follows that 
\begin{align*}
\mathcal{R}_{\varepsilon}f & =F_{f,\theta}(\theta+\varepsilon)\\
 & =\left\Vert f\right\Vert _{p}\left(\frac{\left\vert f\right\vert }{\left\Vert f\right\Vert _{p}}\right)^{h(\theta+\varepsilon)p}\frac{f}{\left\vert f\right\vert }\\
 & =\left(\frac{\left\vert f\right\vert }{\left\Vert f\right\Vert _{p}}\right)^{\varepsilon p(\frac{1}{p_{1}}-\frac{1}{p_{0}})}f.
\end{align*}
 Moreover, 
\[
\mathcal{R}_{\varepsilon}:L^{p}=[L^{p_{0}},L^{p_{1}}]_{\theta}\rightarrow[L^{p_{0}},L^{p_{1}}]_{\theta+\varepsilon}=L^{p_{\varepsilon}},
\]
 where 
\[
\frac{1}{p_{\varepsilon}}=\frac{1-(\theta+\varepsilon)}{p_{0}}+\frac{(\theta+\varepsilon)}{p_{1}}=\frac{1}{p}+\varepsilon(\frac{1}{p_{1}}-\frac{1}{p_{0}}).
\]
 Then, the following commutator theorem, which is essentially the
same as Theorem 4 of \cite[p.~147]{is}, can be seen to also follow
from (\ref{iwan}): \textit{Let $T$ be a bounded operator on $L^{p_{j}}$
for $j=0,1$. Then for each $\varepsilon$ such that $0<\theta+\varepsilon<1$
and for 
\[
[\mathcal{R}_{\varepsilon},T]f=\mathcal{R}_{\varepsilon}(Tf)-T(\mathcal{R}_{\varepsilon}f)=\left(\frac{\left\vert Tf\right\vert }{\left\Vert Tf\right\Vert _{p}}\right)^{\varepsilon p(\frac{1}{p_{1}}-\frac{1}{p_{0}})}(Tf)-T\left(\left(\frac{\left\vert f\right\vert }{\left\Vert f\right\Vert _{p}}\right)^{\varepsilon p(\frac{1}{p_{1}}-\frac{1}{p_{0}})}f\right)
\]
 we have 
\[
\left\Vert [\mathcal{R}_{\varepsilon},T]f\right\Vert _{L^{p_{\varepsilon}}}\le c\left|\varepsilon\right|\left\Vert f\right\Vert _{L^{p}},
\]
 where $c$ depends only on $p$, $p_{0}$, $p_{1}$, $\left\Vert T\right\Vert _{L^{p_{0}}\to L^{p_{0}}}$
and $\left\Vert T\right\Vert _{L^{p_{1}}\to L^{p_{1}}}$.}\end{example}
\begin{rem} \label{rem: ApplThree}Theorem 4 of \cite{is} plays
an important role in that paper. It is a tool there for obtaining
results concerning existence and regularity properties of solutions
(including vector valued ones) of the Dirichlet problem on certain
domains $W\subset\mathbb{R}^{n}$, of a large class of divergence
type partial differential equations. That class includes the familiar
$p$-harmonic equation $\operatorname{div}(\left\vert \nabla u\right\vert ^{p-2}\nabla u)=0$.
In this application the role of $T$ is played by a certain singular
integral operator $\mathcal{K}_{W}$ associated with the so-called
Hodge decomposition of vector fields defined on $W$. On page 147
of \cite{is} there is also a brief indication of other possible significant
applications of Theorem 4 of that paper to nonlinear potential theory
and other topics. Theorem 13.2.1 on p.~341 of the book \cite{IwaniecTMartinG2001}
is also essentially the same as the theorem formulated just above,
and is presented there with rather different applications in mind.
More generally, Chapter 13 of \cite{IwaniecTMartinG2001} offers a
number of applications of the ideas we have been discussing to problems
in analysis. \end{rem}

\subsection{\label{sub:Furthermore}Furthermore}

In \cite{ka} 
{(see also our introduction to \cite{ka}
in \cite{CwikelMMilmanMRochbergR2014})} Nigel defined \textit{a priori}
a class of operators $\Omega$ that form bounded commutators with
bounded operators on a large class of K\"othe spaces. Moreover, in \cite{kami},
he and Marius Mitrea considered a very general interpolation method
for quasi-Banach spaces, which contains both complex interpolation
and real interpolation as special cases, and has an in-built cancellation
principle (expressed by Property (3) on p.~3905 of \cite{kami}).
Their applications included a study of the stability of isomorphism
or Fredholm properties etc.~of operators on complex interpolation
scales of quasi-Banach spaces. By ``stability'' we mean results
of the following kind: Suppose that the spaces $A_{t}$ form an interpolation
scale with respect to the complex (or some other) method, and that
the linear operator $T$ is bounded on each $A_{t}$ and has some
additional property (such as being a Fredholm operator) on $A_{t}$
when $t$ takes some particular value $t=t_{0}$. Then there exists
some open neighbourhood of $t_{0}$ such that $T$ has that same additional
property on $A_{t}$ for all $t$ in that neighbourhood.

\subsection{Methods of real interpolation, and higher order commutators }

The ideas of \cite{rw} were extended to the real method of interpolation
by Jawerth-Rochberg-Weiss \cite{jrw}. In that paper, suitable operators
$\Omega$ were defined for the $K$, $J$ and $E$ methods of interpolation,
and corresponding new commutator theorems were obtained for them.
Explicit calculations of these $\Omega$ operators showed that they
depend on the particular method of interpolation used. It was also
shown that for the $J$-method one can provide a treatment of the
theory of commutators that closely follows the analysis of the complex
method. The form and effects of cancellations for the real method
of interpolation were discussed in some detail in \cite{mr}.

Another natural issue that arose at the time was the problem of dealing
with higher order commutators. Simple iteration does not work since
the (first order) commutator $[T,\Omega]$ is, in general, not linear.
Nevertheless, higher order versions of the commutator theorems for
the complex method were obtained in \cite{r1} and the corresponding
results for the real methods followed suit (cf.~\cite{R5} for a
recent survey).

\section{\label{sec:UnifiedMethod}A unified method of interpolation}

Now we come to the review of the paper \cite{ckmr} proper. As the
interpolation theory of commutators was developing, many of the results
had to be derived, one method of interpolation at a time, and via
different arguments tailored for each specific method of interpolation.
The treatment of the complex method was particularly elegant and amenable
to the study of cancellations, higher order inequalities, and further
applications. However, the study of commutators in the contexts of
other methods of interpolation was somewhat more laborious.

This prompted the following question: Might it be possible to find
some unified way of describing several different interpolation methods,
in a manner which also provides a unified way of obtaining diverse
known (and hopefully also new) commutator results for those different
methods? Our work \cite{ckmr} in 2002 with Nigel%
\footnote{Throughout the rest of this document, please understand ``we'',
``our'' or ``ours'' to often be referring to \textit{all} authors
of \cite{ckmr}, of course including Nigel. %
} provided an affirmative answer to this question (and also showed
that the known straightforward approach for the complex method could
be adapted to work for other methods). But ours was not the first
such answer.

In 1981 the impressive work of Svante Janson \cite{ja} had shown
that most known interpolation methods fit into the general Aronszajn-Gagliardo
framework of orbital and co-orbital methods (see also \cite{ov}).
But at that time there did not seem to be any way of dealing with
commutators in a context as general as this. (See however item 4 in
Section \ref{sec:FurtherDirections} regarding some progress decades
later towards such a goal.)

In 1995 Maria Carro, Joan Cerd\`a and Javier Soria \cite{ccs}, using
a construction due to Vernon Williams \cite{w} as their starting
point, developed a quite general interpolation theory which did enable
treatment of commutators, and in which known commutator theorems,
for example, for the complex method and for the real $J$ and $K$
methods all readily appeared as special cases. They gave further discussion
and development of their work in \cite{barc-2,barc-h,barc-m}. In
several places in our paper, in particular in Section 5 on pp.~276--279,
the reader can find detailed comparisons of their approach and ours
to this topic.

Our paper was dedicated to Jaak Peetre on the occasion of his 65th
birthday, (and the dedication even appears in his native Estonian%
\footnote{It was kindly translated into Estonian for us by the late Simson Baron.%
}). We felt this to be particularly appropriate, since Jaak's 1971
paper \cite{pe} had essentially shown us the way to give several
other interpolation methods a holomorphic ``structure'' similar,
one might even say surprisingly similar, to that which appears intrinsically
in the complex method. Let us indicate how this happened: In \cite{pe},
among many other things, Jaak introduced a new method of interpolation,
sometimes since called the $\pm$ method. The spaces which it generates
are usually denoted by $\left\langle A_{0},A_{1}\right\rangle _{\theta}$.
Jaak observed that they are related to real and complex interpolation
spaces by the inclusions $\left(A_{0},A_{1}\right)_{\theta,1}\subset\left\langle A_{0},A_{1}\right\rangle _{\theta}\subset\left[A_{0},A_{1}\right]_{\theta}$.
In fact in his discussion of these matters he worked with a certain
variant of the complex method, about which we will now need to say
a few words. We will temporarily use the notation $\left[A_{0},A_{1}\right]_{\theta}^{(\mathbb{A})}$
for the spaces obtained by this different method. Let $\mathbb{S}=\left\{ z\in\mathbb{C}:0<\mathrm{Re}\, z<1\right\} $
and $\mathbb{A}=\left\{ z\in\mathbb{C}:1<\left|z\right|<e\right\} $.
It is very easy to see that the definition of the spaces $\left[A_{0},A_{1}\right]_{\theta}^{(\mathbb{A})}$
as formulated in \cite{pe} can be equivalently reformulated to be
the same as Calder\'on's definition in \cite{ca}, except that the closed
strip $\overline{\mathbb{S}}=\left\{ z\in\mathbb{C}:0\le\mathrm{Re}\, z\le1\right\} $,
which was used by Calder\'on in \cite{ca} (and originally by Thorin),
must be replaced in a rather obvious way, by a closed annulus, say%
\footnote{To make our definition here correspond \textit{exactly} to the definition
in \cite{pe} the outer radius of our annulus should be $2$ instead
of $e$. But (see the next footnote) the choice of radius does not
really matter.%
} $\overline{\mathbb{A}}=\left\{ z\in\mathbb{C}:1\le\left|z\right|\le e\right\} $.
More explicitly, let $\mathcal{F}_{\mathbb{A}}(A_{0},A_{1})$ be the
space of continuous functions $f:\overline{\mathbb{A}}\to A_{0}+A_{1}$
which are holomorphic in $\mathbb{A}$ and, for $j=0,1$, their restrictions
to the circles $\left|z\right|=e^{j}$ are continuous $A_{j}$-valued
functions. Let us norm $\mathcal{F}_{\mathbb{A}}(A_{0},A_{1})$ by
\[
\left\Vert f\right\Vert _{\mathcal{F}_{\mathbb{A}}(A_{0},A_{1})}=\max_{j=0,1}\sup\left\{ \left\Vert f(e^{j+it})\right\Vert _{A_{j}}:0\le t\le2\pi\right\} .
\]
 Now we can define $\left[A_{0},A_{1}\right]_{\theta}^{(\mathbb{A})}$
to be the space of all elements $a\in A_{0}+A_{1}$ of the form $a=f\left(e^{\theta}\right)$
for some $f\in\mathcal{F}_{\mathbb{A}}(A_{0},A_{1})$. The norm in
$\left[A_{0},A_{1}\right]_{\theta}^{(\mathbb{A})}$ is defined exactly
analogously to Calder\'on's definition of the norm in $\left[A_{0},A_{1}\right]_{\theta}$,
i.e., by 
\begin{equation}
\left\Vert a\right\Vert _{\left[A_{0},A_{1}\right]_{\theta}^{(\mathbb{A})}}=\inf\left\{ \left\Vert f\right\Vert _{\mathcal{F}_{\mathbb{A}}(A_{0},A_{1})}:f\in\left\Vert f\right\Vert _{\mathcal{F}_{\mathbb{A}}(A_{0},A_{1})},\, a=f(e^{\theta})\right\} .\label{eq:AnnNorm}
\end{equation}
 In \cite{pe} and elsewhere, Jaak Peetre asked whether $\left[A_{0},A_{1}\right]_{\theta}^{(\mathbb{A})}$
is the same space as $\left[A_{0},A_{1}\right]_{\theta}$. The answer
\cite{cw} turned out to be yes, they are equal%
\footnote{The same result holds when the annulus defining this space is replaced
by $\overline{A(r)}:=\left\{ z:1\le\left|z\right|\le r\right\} $
for any choice of $r>1$. The norm of the complex interpolation space
defined using $\overline{\mathbb{A}(r)}$ is only known to be equivalent
to the norm of the space $\left[A_{0},A_{1}\right]_{\theta}$, and
it is known that in general these norms cannot be equal, e.g.~for
large values of $r$. Recently Eliran Avni \cite{AvniE2014} has shown
that the norm of this space becomes arbitrarily close to the norm
of $\left[A_{0},A_{1}\right]_{\theta}$ as $r$ tends to $1$. %
}, to within equivalence of norms.

Now let us rewrite the above definition of $\left[A_{0},A_{1}\right]_{\theta}^{(\mathbb{A})}$
in terms of the coefficients which appear in the Laurent expansions
$f(z)=\sum_{n\in\mathbb{Z}}z^{n}a_{n}$ in $\overline{\mathbb{A}}$
of functions $f\in\mathcal{F}_{\mathbb{A}}(A_{0},A_{1})$. It is easy
to show that the coefficients $a_{n}$ in such expansions must all
be elements of $A_{0}\cap A_{1}$.{} Our equivalent definition has
to be as follows: 
\begin{defn}
\label{def:ViaLaurent}The space $\left[A_{0},A_{1}\right]_{\theta}^{(\mathbb{A})}$
consists of all elements $a\in A_{0}+A_{1}$ which can be written
in the form 
\begin{equation}
a=\sum_{n\in\mathbb{Z}}e^{\theta n}a_{n}\label{eq:asum}
\end{equation}
 where $\left\{ a_{n}\right\} _{n\in\mathbb{Z}}$ is a sequence of
elements in $A_{0}\cap A_{1}$ such that, for $j=0$ and $j=1$, the
``weighted'' sequence $\left\{ e^{jn}a_{n}\right\} _{n\in\mathbb{Z}}$
belongs to the sequence space $FC\left(A_{j}\right)$. 
\end{defn}
Here, for any Banach space%
\footnote{i.e., any \textit{complex} Banach space. All Banach spaces in this
discussion are assumed to be complex.%
} $B$, the sequence space $FC(B)$ is defined to consist of all $B$-valued
sequences $\left\{ b_{n}\right\} _{n\in\mathbb{Z}}$ which arise as
the Fourier coefficients of a (uniquely determined) continuous $B$-valued
function $\phi:\mathbb{T}\to B$, and the norm $\left\Vert \left\{ b_{n}\right\} _{n\in\mathbb{Z}}\right\Vert _{FC(B)}$
of such a sequence equals $\sup\left\{ \left\Vert \phi(e^{it})\right\Vert _{B}:t\in[0,2\pi]\right\} $.

It is convenient to let $\mathcal{J}(FC(A_{0}),FC(A_{1}))$ denote
the space of all sequences $\left\{ a_{n}\right\} _{n\in\mathbb{Z}}$
which satisfy the two conditions specified in Definition \ref{def:ViaLaurent}
and to equip it with the norm 
\begin{equation}
\left\Vert \left\{ a_{n}\right\} _{n\in\mathbb{Z}}\right\Vert _{\mathcal{J}(FC(A_{0}),FC(A_{1}))}=\max_{j=0,1}\left\Vert \left\{ e^{jn}a_{n}\right\} _{n\in\mathbb{Z}}\right\Vert _{FC(A_{j})}.\label{eq:djnorm}
\end{equation}
 So $\mathcal{J}(FC(A_{0}),FC(A_{1}))$ is simply the ``Fourier transform''
of $\mathcal{F}_{\mathbb{A}}(A_{0},A_{1})$ and we can rewrite (\ref{eq:AnnNorm})
as 
\begin{equation}
\left\Vert a\right\Vert _{\left[A_{0},A_{1}\right]_{\theta}^{(\mathbb{A})}}=\inf\left\{ \left\Vert \left\{ a_{n}\right\} _{n\in\mathbb{Z}}\right\Vert _{\mathcal{J}(FC(A_{0}),FC(A_{1}))}\right\} \label{eq:jnorm}
\end{equation}
 where the infimum is taken over all sequences $\left\{ a_{n}\right\} _{n\in\mathbb{Z}}$
in $\mathcal{J}(FC(A_{0}),FC(A_{1}))$ which satisfy (\ref{eq:asum}).

The definition of $\left[A_{0},A_{1}\right]_{\theta}^{(\mathbb{A})}$,
when expressed in the language of Definition \ref{def:ViaLaurent},
looks very much like one of the equivalent definitions of the Lions-Peetre
real interpolation space which is usually denoted by $\left(A_{0},A_{1}\right)_{\theta,p}$.
(Cf.~e.g.~\cite[Lemma 3.2.3  p.~43]{bl} or \cite[p.~17]{lp}, where
different notation is used for this space.) If we simply replace the
spaces $FC(A_{0})$ and $FC(A_{1})$ respectively by $\ell^{p}(A_{0})$
and $\ell^{p}(A_{1})$ then it becomes exactly that definition. (Here,
as usual, $\ell^{p}(B)$ is the space of all $B$-valued sequences
$\left\{ b_{n}\right\} _{n\in\mathbb{Z}}$ for which $\left(\sum_{n\in\mathbb{Z}}\left\Vert b_{n}\right\Vert _{B}^{p}\right)^{1/p}$
is finite.) Furthermore, replacing $FC(A_{0})$ and $FC(A_{1})$ by
$\ell^{p}(A_{0})$ and $\ell^{p}(A_{1})$ in the definition of $\mathcal{J}(FC(A_{0}),FC(A_{1}))$
and then also in (\ref{eq:jnorm}), converts (\ref{eq:jnorm}) (of
course together with the words immediately following it) into a formula
for one of the equivalent norms of $\left(A_{0},A_{1}\right)_{\theta,p}$.

Yet another analogous change will convert Definition \ref{def:ViaLaurent}
and the accompanying definition of the norm of $\left[A_{0},A_{1}\right]_{\theta}^{(\mathbb{A})}$,
into the definition of the space $\left\langle A_{0},A_{1}\right\rangle _{\theta}$
and its norm. This time $FC(A_{0})$ and $FC(A_{1})$ must be replaced
by the sequence spaces $UC(A_{0})$ and $UC(A_{1})$ whose definitions
can be found in \cite[Example 2.4, p.~247]{ckmr}. The Gustavsson-Peetre
variant of $\left\langle A_{0},A_{1}\right\rangle _{\theta}$ (introduced
in \cite{gp}) will be obtained if we replace $FC(A_{0})$ and $FC(A_{1})$
by certain other sequence spaces (\cite[Example 2.5, p.~247]{ckmr}).

Almost all of the above was explicit or implicit in \cite{pe}. The
first step that we took in \cite{ckmr} (in Section 2 of that paper)
was to incorporate all of the interpolation spaces given by the above
definitions, and also other more exotic ones, into the general framework
of \textit{interpolation spaces defined via pseudolattices}. We defined
a pseudolattice to be a map (or one might say functor) $\mathcal{X}$
which, when applied to any complex Banach space $B$, yields a certain
sequence space $\mathcal{X}(B)$ of $B$-valued sequences. This map
was required to satisfy some simple conditions which we will not list
here. (See \cite[Definition 2.1 on p.~246]{ckmr}.) For example, if
we choose some suitable Banach lattice $X$ on $\mathbb{Z}$, whose
elements can thus be considered as complex valued sequences $\left\{ \alpha_{n}\right\} _{n\in\mathbb{Z}}$,
then we can use $X$ to define a particular pseudolattice $\mathcal{X}$
by requiring that, for each Banach space $B$, the space $\mathcal{X}(B)$
is the space (often denoted by $X(B)$) of all sequences $\left\{ b_{n}\right\} _{n\in\mathbb{Z}}$
in $B$ for which $\left\{ \left\Vert b_{n}\right\Vert _{B}\right\} _{n\in\mathbb{Z}}\in X$.
This is exactly what we did above in the special case $X=\ell^{p}$.
In fact our discussion just above implicitly defines three pseudolattices:
$\ell^{p}$ and also $FC$ and $UC$.

Our definition of pseudolattice was formulated so that each choice
of pseudolattice $\mathcal{X}$ and parameter $\theta\in(0,1)$ would
provide us with a new interpolation method. This method, when applied
to any Banach pair $\left(A_{0},A_{1}\right)$ yields an interpolation
space which will be denoted here%
\footnote{In order to simplify our presentation here we are using notation which
is a little different from that of \cite{ckmr} and not considering
the more general version of this definition there. In that version
we had the additional options of applying \textit{different} pseudolattices
to $A_{0}$ and $A_{1}$ and of permitting the parameter $\theta$
to also have a non-zero imaginary part. I.e., there we replaced the
parameter $\theta\in(0,1)$ by the parameter $s=e^{\theta}\in\mathbb{A}$.%
}~by $\left(A_{0},A_{1}\right)_{\mathcal{X},\theta}$. As can be anticipated
from the preceding discussion, its definition and its norm were obtained
by simply replacing $FC(A_{0})$ and $FC(A_{1})$ by $\mathcal{X}(A_{0})$
and $\mathcal{X}(A_{1})$ throughout the preceding definitions.

In order to get to this stage, it had been convenient to make the
transition from holomorphic functions $f:\mathbb{A}\to A_{0}+A_{1}$,
to $A_{0}\cap A_{1}$ valued sequences $\left\{ a_{n}\right\} _{n\in\mathbb{Z}}$
(the Laurent coefficients of $f$). But then we wanted to go from
these sequences back to holomorphic functions on $\mathbb{A}$. It
turned out that a certain mild condition (which we called \textit{Laurent
compatibility}%
\footnote{The original definition of Laurent compatibility in \cite{ckmr} (see
Definition 2.9 on p.~248) was formulated for \textit{pairs} of pseudolattices.%
}) on $\mathcal{X}$, (and which holds in all ``natural'' examples
known so far) suffices to ensure that, for every sequence $\left\{ a_{n}\right\} _{n\in\mathbb{Z}}$
in $\mathcal{J}(\mathcal{X}(A_{0}),\mathcal{X}(A_{1}))$, the series
$\sum_{n\in\mathbb{Z}}z^{n}a_{n}$ converges for all $z\in\mathbb{A}$
to a holomorphic $A_{0}+A_{1}$-valued function.

Now things began to look good. (We are now talking about Section 3
of \cite{ckmr}.) For each $\theta\in(0,1)$ we expected to be able
to define an operator $\Omega$ for $A_{\theta}=\left(A_{0},A_{1}\right)_{\mathcal{X},\theta}$
very similarly to how this had been done in (\ref{eq:DefOmega}) and
then to prove a more general result like (\ref{cuatro}) by the same
kind of simple argument that was sketched above just after (\ref{cuatro}).
We also expected to be able to define an analogue of the operator
$\mathcal{R}_{\varepsilon}$ of Subsection \ref{sub:transops} and
obtain a result like (\ref{iwan}) for it. Indeed (via our statement
and proof of Theorem 3.8 of \cite[p.~256]{ckmr}), we would accomplish
both these tasks. But there were two issues to be dealt with before
this could happen.

The first of these was that (as indeed in the previously considered
special cases of such results) we could not be sure that in general
the infimum in the generalized version of (\ref{eq:jnorm}) is attained
by some sequence. (In fact an analogous concern arises already for
the prototype definition in (\ref{eq:DefOmega}).) Even if this infimum
is attained, the sequence attaining it may fail to be unique. So how
can we define $\Omega a$ or $\mathcal{R}_{\varepsilon}a$ for each
$a\in\left(A_{0},A_{1}\right)_{\mathcal{X},\theta}$? We overcame
or bypassed this difficulty rather easily, by allowing $\Omega a$
and $\mathcal{R}_{\varepsilon}a$ to be appropriate \textit{sets}
of elements, rather than single elements. As stated in \cite[Definition 3.1, p.~253]{ckmr},
for some chosen constant $C_{opt}>1$, these sets are, respectively,
the set of all values of $f^{\prime}(e^{\theta})$ or of $f\left(e^{\theta+\varepsilon}\right)$
as $f$ ranges over all holomorphic functions $f:\mathbb{A}\to A_{0}+A_{1}$
which satisfy $f\left(e^{\theta}\right)=a$, and whose sequences of
Laurent coefficients have norms in $\mathcal{J}(\mathcal{X}(A_{0}),\mathcal{X}(A_{1}))$
which are no greater than $C_{opt}$ times the infimum defining the
norm $\left\Vert a\right\Vert _{\left(A_{0},A_{1}\right)_{\mathcal{X},\theta}}$
of $a$. (Later below there will be some brief discussion of $\Omega_{n}$,
an ``$n$-order'' version of $\Omega$. For each $a$ as above,
$\Omega_{n}a$ is obtained by replacing $f^{\prime}(e^{\theta})$
by $f^{(n)}(e^{\theta})$ in the above definition.)

To explain the second issue, one should look again at the above-mentioned
proof of the special case sketched just after (\ref{cuatro}). Suppose
that $g(z)=\sum_{n\in\mathbb{Z}}z{}^{n}a_{n}$ for some sequence $\left\{ a_{n}\right\} _{n\in\mathbb{Z}}$
in $\mathcal{J}(\mathcal{X}(A_{0}),\mathcal{X}(A_{1}))$ and that
$g(e^{\theta})=0$. Here, analogously to the last step (immediately
after (\ref{eq:vjta})) of` that proof, we needed to know that, if
we multiply $g(z)$ by the scalar function $1/\left(z-e^{\theta}\right)$,
then we will obtain a new holomorphic function on $\mathbb{A}$ whose
Laurent coefficients keep the same property as that possessed by those
of $g$, i.e., they too are a sequence in $\mathcal{J}(\mathcal{X}(A_{0}),\mathcal{X}(A_{1}))$.
It turned out that there is a very simple condition which suffices
to ensure that this happens, and this condition is satisfied when
$\mathcal{X}$ is $FC$ or $\ell^{p}$ or $UC$ or any one of a large
number of other ``natural'' choices. I.e., it suffices if, for each
Banach space $B$, the left shift operator $\left\{ b_{n}\right\} _{n\in\mathbb{N}}\mapsto\left\{ b_{n+1}\right\} _{n\in\mathbb{N}}$
maps the sequence space $\mathcal{X}(B)$ isometrically onto itself.
(Then of course the right shift operator $\left\{ b_{n}\right\} _{n\in\mathbb{N}}\mapsto\left\{ b_{n-1}\right\} _{n\in\mathbb{N}}$
has the same property.)

In fact we could also manage without imposing this isometry of shifts
condition on $\mathcal{X}$. It turned out that a much weaker but
more technical condition%
\footnote{We used the terminology \textit{$\mathcal{X}$ admits differentiation}
for this property. In fact (see Definition 3.4 of \cite[pp.~254--255]{ckmr})
we defined and dealt with \textit{pairs} of pseudolattices which admit
differentiation. Cf.~also Lemma 3.6 of \cite{ckmr} p.~255 for another
condition which can work in place of isometry of shifts.%
} on $\mathcal{X}$ is also sufficient to ensure that multiplication
of $g(z)$ by $1/(z-e^{\theta})$ has the property that we require.

Now, having described the path towards it, and the conditions required
to take that path, we can invite you to look at Theorem 3.8 of \cite[p.~256]{ckmr}.
As you can see, it includes results mentioned above as special cases.
In particular, the inequalities (3.2) and (3.3) in the statement of
that theorem indeed generalize (\ref{cuatro}) and (\ref{iwan}),
respectively. Furthermore, as we showed after proving it, (see Section
4 of \cite{ckmr}) a number of other previously known results for
commutators, including some in the context of versions of the real
interpolation method, can also be deduced from Theorem 3.8. \begin{rem}
\label{rem:ApplTwo}The applicability of special cases of Theorem
3.8 to various topics in analysis (cf.~Remark \ref{rem:ApplOne}
and the comments near the beginning of Subsection \ref{sub:transops})
suggest that the similar but more general perspectives offered by
Theorem 3.8 should also have useful applications. \end{rem} It remains
to give some brief account of what can be found in the several subsequent
sections of \cite{ckmr}. As already mentioned, Section 5 compared
our approach with that of Carro-Cerd\`a-Soria. Then in Section 6 we
turned to a study of higher order analogs of Theorem 3.8. The natural
definition of $\Omega_{n}$, the $n$th order analog of $\Omega$
has already been mentioned above. For bounded operators $T$, our
first theorem in Section 6 did not control the norm of $\left[T,\Omega_{n}\right]a$
but rather of a more elaborate inductively defined expression which
contains ``correction terms'' also involving $\Omega_{k}$ for $k=1,2,...,n-1$.
Analogs of this result had been obtained earlier in \cite{barc-2}
and also, in the particular contexts of real or complex interpolation
methods, in \cite{r1} and \cite{mi}. (Cf.~also \cite{ccms} and
\cite{mr}.) We also obtained a second theorem, this time for translation
operators $\mathcal{R}_{\varepsilon}$, which generalized the pseudolattice
version of the norm estimate (\ref{iwan}). Here the factor $\left|\varepsilon\right|$
on the right side of the estimate could be replaced by $\left|\varepsilon\right|^{n}$.
Analogously to the preceding theorem for $\Omega_{n}$, in this theorem
the commutator expression of $T$ and $\mathcal{R}_{\varepsilon}$
on the left side whose norm is controlled by this power of $\left|\varepsilon\right|$
here also has to be quite elaborate and include extra ``correction
terms'' which also involve $\Omega_{k}$ for $k=1,2,..,n-1$.

As indicated above (in particular in Remarks \ref{rem:ApplOne} and
\ref{rem: ApplThree}) nonlinear first order commutators have already
found quite a range of applications. It seems reasonable to expect
that higher order results, such as those that we have mentioned here,
will also eventually find interesting applications in various branches
of analysis. For example, perhaps they could turn out to be useful
in the study of functional equations (cf.~\cite{KoMi} and the references
therein).

In Section 7 of \cite{ckmr} we extended results about domain spaces
and range spaces of operators $\Omega$ previously known for the real
and complex methods to the general pseudolattice method. In both this
and the previous section we noted connections with the Lions-Schechter
variant \cite{schechterm} of the complex interpolation spaces. (We
showed later, in an appendix, that also for the definition of these
spaces, one can replace the strip by an annulus.)

Finally, to describe the result of Section 8, we first need to recall
that there are two main ways of constructing the Lions-Peetre real
interpolation spaces $\left(A_{0},A_{1}\right)_{\theta,p}$. These
are usually referred to as the $J$-method and the $K$-method. The
fact that these two methods give the same spaces is sometimes quite
helpful, for example for characterizing the duals of these spaces
or for proving reiteration formulae (e.g.~for describing the space
$\left(\left(A_{0},A_{1}\right)_{\theta_{0},p_{0}},\left(A_{0},A_{1}\right)_{\theta_{1},p_{1}}\right)_{\theta_{3},p_{3}}$).
The construction of the space $\left(A_{0},A_{1}\right)_{\mathcal{X},\theta}$
for an arbitrary pseudolattice $\mathcal{X}$ is modeled on the $J$-method
(which explains our choice of notation for the related space $\mathcal{J}(\mathcal{X}(A_{0}),\mathcal{X}(A_{1}))$.)
Somewhat to our surprise, it turned out to be possible to formulate
an alternative way of constructing the general spaces $\left(A_{0},A_{1}\right)_{\mathcal{X},\theta}$
which is modeled on the $K$-method for $\left(A_{0},A_{1}\right)_{\theta,p}$.
We are not aware of this fact having been observed earlier, not even
for the important and much studied special case of complex interpolation
spaces $\left[A_{0},A_{1}\right]_{\theta}$. Results about duality
and reiteration are well known for $\left[A_{0},A_{1}\right]_{\theta}$.
But we keep the hope that this unexpected $K$-method will prove to
be useful for other purposes.

\section{\label{sec:FurtherDirections}Further directions for future research}

In this final section we offer some additional comments, and formulations
of some natural open problems associated with results discussed in
the preceding sections. 
\begin{itemize}
\item [(1)] The first project asks for a version of the interpolation
of analytic families of operators of Stein \cite{SteinE1956} (cf.~also
\cite{cj}), in the context of the interpolation spaces of \cite{ckmr}.
(Cf.~also the notion of ``interpolating family of operators'' introduced
in \cite[p.~3905]{kami}.) 
\item [(2)] Other open problems for the interpolation method of \cite{ckmr}
include: the duality theory, bilinear interpolation and the study
of compactness for these methods. When dealing with these it may turn
out to be helpful to know that a lemma of Stafney, which is sometimes
useful for the study of the complex method has been shown \cite{Ivtsan}
to also hold for this general method. 
\item [(3)] Much remains to be done in the theory of interpolation of
commutators. In particular, the problem of compactness of commutators
has not been treated in the context of the abstract theory. The classical
result due to Uchiyama \cite{UchiyamaA1978} which one could seek
to generalize states that if $T$ is a Calderón-Zygmund operator and
$b\in VMO,$ and $p\in(1,\infty),$ then the commutator $[T,M_{b}]:L^{p}\rightarrow L^{p}$
is compact. What are the compactness results in the general theory?
The problem is wide open even for the classical methods of interpolation.
(Note that B\'enyi and Torres \cite{BenyiATorresR2013} have generalized
Uchiyama's result \cite{UchiyamaA1978} to the context of bilinear
operators.) 
\item [(4)] We have seen that cancellation principles can be formulated
and used to obtain commutator theorems in quite general contexts,
as in \cite{ccs}, \cite{ckmr}, \cite{kami}, etc. But still greater
generality is possible. An abstract cancellation principle for orbital
interpolation methods generated by a single element appears in \cite{krmi}
and is applied, not only to commutators, but also to stability%
\footnote{Here ``stability'' has the same meaning as specified in Subsection
\ref{sub:Furthermore}.%
} theorems. Such orbital methods can be seen to extend some of the
classical methods. (See \cite{ja,krmi}). But are there corresponding
results for orbital methods with a finite number of generators or,
better still, without any such restrictions? Likewise what is the
corresponding theory of commutators for co-orbital methods? (The work
of \cite{krmi} has a sequel \cite{KaspikrovaE2009} which presents
some higher order results.) 
\item [(5)] Establish bilinear commutator theorems. The problem seems
to be open even for the classical methods of interpolation. As already
mentioned in item (3), some specialized bilinear results can be found
in the paper \cite{BenyiATorresR2013}. 
\item [(6)] The commutator theorem of \cite{crw} involves the commutator
of singular integrals and the multiplication operator $M_{b},$ with
$b\in BMO$. The method of \cite{rw} also gives the $L^{p}$ boundedness
of commutators $[T,M_{b}]$ for all operators $T$ that are bounded
on $L^{p}(w)$ for all weights in the Muckenhoupt class $A_{p}$%
\footnote{This follows by the same reasoning as that used in Example \ref{ex:bmo}
above.%
}. More generally, the result can be extended to some nonlinear operators
$T$ that satisfy the same weighted norm inequalities. For example
(cf.~\cite{bmrproc} and the references therein), if $T$ is either
the maximal operator of Hardy-Littlewood or the sharp maximal operator,
then $[T,M_{b}]$ is bounded on $L^{p}$ if and only if $b\in BMO$
and its negative part $b^{-}$ is bounded%
\footnote{In particular if $b\in BMO$ and $b\geq0.$%
}$.$ While some other results are known for quasilinear operators
(cf.~\cite{cjmr}) it would be of interest to have a fully developed
systematic theory of nonlinear commutators. 
\item [(7)] The problems that we propose here are motivated by an interesting
application of commutator theorems given in \cite{Bon} in connection
with the definition of the product operator $M_{b}$ in $H^{1}$.
One learns of the important roles in concrete applications of Orlicz
spaces as well as Lorentz and Marcinkiewicz spaces in the commutator
theory of interpolation. These give us good reasons to consider the
project of finding an extension of the theory of \cite{ckmr} to be
able to incorporate such refinements. 
\end{itemize}
\medskip{}

\textbf{\textit{Acknowledgment:}} We thank Pavel Shvartsman for some
helpful correspondence.

\bigskip{}


\begin{thebibliography}{10}
\bibitem{AvniE2014}E.~Avni, The periodic complex method in interpolation
spaces. \textit{Rev.~Mat.~Complut.} \textbf{27} (2014), 156--166.
Earlier version: \texttt{arXiv:1203.4183 {[}math.FA{]}}

\bibitem{bmrproc}J.~Bastero, M.~Milman and F.~J.~Ruiz, Commutators
for the maximal and sharp functions, \textit{Proc.~Amer.~Math.~Soc.}
\textbf{128} (2000), 3329--3334.

\bibitem{BenyiATorresR2013}A.~B\'enyi and R.~H.~Torres, Compact
bilinear operators and commutators, \textit{Proc.~Amer.~Math.~Soc.}
\textbf{141} (2013), 3609--3621.

\bibitem{bl}J.\ Bergh and J.\ Löfström, \textit{Interpolation spaces.
An Introduction,} Springer, Berlin 1976.

\bibitem{Bon}A.~Bonami, T.~Iwaniec, P.~Jones, and M.~Zinsmeister,
On the product of functions in BMO and $H^{1}$, \textit{Ann.~Inst.~Fourier}
\textbf{57} (2007), 1405--1439.

\bibitem{ca}A.\ P.\ Calderón, Intermediate spaces and interpolation,
the complex method, \textit{Studia Math.}\ \textbf{24} (1964) 113--190.

\bibitem{ccms}M.\ Carro, J.\ Cerdà, M.\ Milman and J.\ Soria,
Schechter methods of interpolation and commutators, \textit{Math.~Nachr.}
\textbf{174} (1995), 35--53.

\bibitem{ccs}M.~J.~Carro, J.~Cerdà and J.~Soria, Commutators
and interpolation methods, \textit{Ark. Mat.} \textbf{33} (1995),
199--216.

\bibitem{barc-2}M.\ Carro, J.\ Cerdà and J.\ Soria, Higher order
commutators in interpolation theory, \textit{Math.\ Scand.} \textbf{77}
(1995) 301--309.

\bibitem{barc-m}\ M.\ Carro, J.\ Cerdà and J.\ Soria, A unified
approach for commutator theorems in interpolation theory, a survey.
\textit{Rev.\ Mat.\ Univ.\ Complut. Madrid} \textbf{9} (1996),
Special Issue. suppl., 91--108.

\bibitem{barc-h}\ M.\ Carro, J.\ Cerdà and J.\ Soria, Commutators,
interpolation and vector function spaces, Function spaces, interpolation
spaces, and related topics (Haifa 1995), 24--31, \textit{Israel Math.\ Conf.\ Proc.,}
\textbf{13}, Bar-Ilan Univ., Ramat Gan, 1999.

\bibitem{CoiFef}R.~R.~Coifman and C.~Fefferman, Weighted norm
inequalities for maximal functions and singular integrals, \textit{Studia
Math.} \textbf{51} (1974), 241--250.

\bibitem{CLMS}R.~R.~Coifman, P.-L.~Lions, Y.~Meyer, and S.~Semmes,
Compensated compactness and Hardy spaces. \textit{J.~Math.~Pures
Appl.} (9) \textbf{72} (1993), no. 3, 247--286. %46E99 (35S05 42B30 46F10 46N20 47G30)



\bibitem{crw}R.\ R.\ Coifman, R.\ Rochberg and G.\ Weiss, Factorization
theorems for Hardy spaces in several variables, \textit{Ann.\ Math.}
\textbf{103} (1976), 611--635.

\bibitem{cw}M.\ Cwikel, Complex interpolation, a discrete definition
and reiteration, \textit{Indiana Univ.\ Math.\ J.\ }\textbf{27}
(1978), 1005--1009.

\bibitem{cjmr}M.\ Cwikel, B.\ Jawerth, M.\ Milman and R.\ Rochberg,
Differential estimates and commutators in interpolation theory, ``Analysis
at Urbana II'', \textit{London Math.\ Soc., Lecture Note Series,}
vol.\ 138, Cambridge Univ.\ Press, 1989, pp.\ 170--220.

\bibitem{cj}M.~Cwikel and S.~Janson, Interpolation of analytic
families of operators, \textit{Studia Math.} \textbf{79} (1984), 61--71.

\bibitem{ckmr}M.~Cwikel, N.~Kalton, M.~Milman and R.~Rochberg,
A unified theory of commutator estimates for a class of interpolation
methods. \textit{Adv.~Math.} \textbf{169} (2002), no. 2, 241--312.

\bibitem{CwikelMMilmanMRochbergR2014}{M.~Cwikel,
M.~Milman and R.~Rochberg, An introduction to Nigel Kalton's work
on differentials of complex interpolation processes for Kothe spaces.
}\texttt{{arXiv:1404.2893 {[}math.FA{]} }}

\bibitem{MingFAN} M.~Fan, Commutator estimates for interpolation
scales with holomorphic structure. \textit{Complex Anal.~Oper.~
Theory} \textbf{4} (2010), no. 2, 159--178.

\bibitem{GarsiaRubio} J.~García-Cuerva and J.~L.~Rubio de Francia,
\textit{Weighted norm inequalities and related topics.} North-Holland
Mathematics Studies, 116. Notas de Matemática {[}Mathematical Notes{]},
104. North-Holland Publishing Co., Amsterdam, 1985. x+604 pp. ISBN:
0-444-87804-1

\bibitem{gp}\ J.\ Gustavsson and J.\ Peetre, Interpolation of
Orlicz spaces, \textit{Studia Math.} \textbf{60} (1977), 33--59.

\bibitem{HuntMuckenhouptWheeden}R.~Hunt, B.~Muckenhoupt and R.~Wheeden,
Weighted norm inequalities for the conjugate function and Hilbert
transform. \textit{Trans.~Amer.~Math.~Soc.} \textbf{176} (1973),
227--251.

\bibitem{Ivtsan} A.~Ivtsan, Stafney's lemma holds for several ``classical''
interpolation methods. \textit{Proc.~Amer.~Math.~Soc.} \textbf{140}
(2012), no. 3, 881--889.

\bibitem{i0}T. Iwaniec,\textsl{ }$p$-harmonic tensors and quasiregular
mappings, \textit{Ann.~Math.} \textbf{136} (1992), 589-624.

\bibitem{i}T.\ Iwaniec, \textsl{Nonlinear differential forms}, Lectures
in Jyväskylä, Dep.~of Math., Report \textbf{80}, Univ. of Jyväskylä,
Jyväskylä, 1998.

\bibitem{IwaniecTMartinG2001}T.~Iwaniec and G.~Martin, \textit{Geometric
function theory and non-linear analysis.} Oxford Mathematical Monographs.
The Clarendon Press, Oxford University Press, New York, 2001. % xvi+552. ISBN 0-19-850929-4 



\bibitem{is} T.\ Iwaniec and C.\ Sbordone, Weak minima of variational
integrals, \textit{J.\ Reine Angew.\ Math.} \textbf{454} (1994),
143--161.

\bibitem{ja}\ S.\ Janson, Minimal and maximal methods of interpolation,
\textit{J.\ Functional Analysis} \textbf{44} (1981), 50--73.

\bibitem{jrw}B.~Jawerth, R.~Rochberg and G.~Weiss, Commutator
and other second order estimates in real interpolation theory. \textit{Ark.~Mat.}
\textbf{24} (1986), no. 2, 191--219.

\bibitem{ka}N.\ J.\ Kalton, Differentials of complex interpolation
processes for Köthe function spaces, \textit{Trans.\ Amer.\ Math.\ Soc.}
\textbf{333} (1992), 479--529.

\bibitem{kami}N.~Kalton and M.~Mitrea, Stability results on interpolation
scales of quasi-Banach spaces and applications, \textit{Trans.~Amer.~Math.~Soc.}
\textbf{350} (1998), 3903-3922.

\bibitem{KaspikrovaE2009}E.~Kasprikova, Higher order commutator
theorems in the method of orbits, Ph.D. Thesis, Florida Atlantic University,
2009.

\bibitem{KoMi}H.~K\"onig and V.~Milman, Operator equations and domain
dependence, the case of the Schwarzian derivative. \textit{J.~Functional
Analysis} \textbf{266} (2014), 2546--2569.

\bibitem{krmi}N.~Krugljak and M.~Milman, A distance between orbits
that controls commutator estimates and invertibility of operators,
\textit{Adv.~Math.} \textbf{182} (2004), 78--123.

\bibitem{lp}J.\ L.\ Lions and J.\ Peetre, Sur une classe d'espaces
d'interpolation. \textit{Inst.\ Hautes Etudes Sci.\ Publ.\ Math.}
\textbf{19} (1964), 5--68.

\bibitem{mi}M.\ Milman, Higher order commutators in the real method
of interpolation, \textit{Journal D'Analyse Math.} \textbf{37} (1995),
37--55.

\bibitem{mr}M.\ Milman and R.\ Rochberg, The role of cancellations
in interpolation theory, \textit{Contemporary Math.} \textbf{189}
(1995), 403--419.

\bibitem{ov}V.\ I.\ Ovchinnikov, \textit{The method of orbits in
interpolation theory,} Mathematical Reports, Vol.~1, Part 2, Harwood
Academic Publishers 1984, 349--516.

\bibitem{pe}\ J.\ Peetre, Sur l'utilization des suites inconditionellement
sommables dans la th\'eorie des espaces d'interpolation.\ \textit{Rend.\ Sem.\ Mat.\ Univ.\ Padova}
\textbf{46} (1971), 173--190.

\bibitem{rw} R.\ Rochberg and G.\ Weiss, Derivatives of analytic
families of Banach spaces, \textit{Ann.\ Math.} \textbf{118} (1983),
315--347.

\bibitem{r1}R.\ Rochberg, Higher order estimates in complex interpolation
theory, \textit{Pacific J.\ Math.} \textbf{174} (1996), 247--267.

\bibitem{R5}R.~Rochberg, Uses of commutator theorems in analysis.
\textit{Interpolation theory and applications,} 277--295, Contemp.~Math.,
\textbf{445}, Amer.~Math.~Soc., Providence, RI, 2007.

\bibitem{schechterm}M.\ Schechter, Complex interpolation, \textit{Comp.\ Math.}
\textbf{18} (1967), 117--147.\textbf{\textit{ }}

\bibitem{SteinE1956}E.~M.~Stein, Interpolation of linear operators.
\textit{Trans.~Amer.~Math.~Soc.} \textbf{83} (1956), 482--492.

\bibitem{UchiyamaA1978}A.~Uchiyama, On the compactness of operators
of Hankel type. \textit{Tôhoku Math. J. (2)} \textbf{30} (1978), 1630--171.

\bibitem{w}V.\ Williams, Generalized interpolation spaces, \textit{Trans.\ Amer.\ Math.\ Soc.}
\textbf{156} (1971), 309--334. \end{thebibliography}
\end{document}